\documentclass[12pt,a4paper,oneside]{article}

\usepackage{amsthm}
\usepackage{amssymb}
\usepackage{latexsym}
\usepackage{amsmath}
\usepackage{amsfonts}
\usepackage{graphicx}

\usepackage{ifthen}

\theoremstyle{definition}
\newtheorem{defi}{Definition}[section]
\newtheorem{rem}[defi]{\small Remark}
\newtheorem{exam}[defi]{\small Example}
\theoremstyle{plain}
\newtheorem{theo}[defi]{Theorem}
\newtheorem{lem}[defi]{Lemma}
\newtheorem{prop}[defi]{Proposition}
\newtheorem{cor}[defi]{Corollary}

\newcommand{\NN}{\mathbb N}
\newcommand{\ZZ}{\mathbb Z}
\newcommand{\QQ}{\mathbb Q}
\newcommand{\RR}{\mathbb R}
\newcommand{\CC}{\mathbb C}

\newcommand{\PP}{\mathbb P}

\newcommand{\alg}   {\operatorname{alg}}
\newcommand{\aut}   {\operatorname{Aut}}
\newcommand{\degr}  {\operatorname{deg}}
\newcommand{\Dim}   {\operatorname{dim}}

\newcommand{\ext}   {\operatorname{Ext}}
\newcommand{\Hom}   {\operatorname{Hom}}
\newcommand{\norm}     {\operatorname{N}}
\newcommand{\Char}  {\operatorname{char}}
\newcommand{\tor}     {\operatorname{t}}
\newcommand{\Tor}     {\operatorname{T}}

\newcommand{\gal}     {\operatorname{G}}
\newcommand{\ord}     {\operatorname{ord}}
\newcommand{\sep}     {\operatorname{sep}}
\newcommand{\supp}    {\operatorname{supp}}
\newcommand{\I}       {\operatorname{I}}
\newcommand{\ab}      {\operatorname{ab}}
\newcommand{\solv}    {\operatorname{solv}}
\newcommand{\real}    {\operatorname{real}}
\newcommand{\ered}    {\operatorname{er}}
\newcommand{\red}      {\operatorname{r}}
\newcommand{\Ho}       {\operatorname{H}}
\newcommand{\id}      {\operatorname{id}}
\newcommand{\coker}   {\operatorname{coker}}
\newcommand{\her}     {\operatorname{h}}

\newenvironment{myabstract}{\centerline{\large\bfseries Abstract} \medskip
\small}

\newenvironment{myrem}{\begin{rem}\small}{\end{rem}}
\newenvironment{myexam}{\begin{exam}\small}{\end{exam}}

\usepackage{titlesec}
\titleformat{\section}%
  {\normalfont \bfseries \large}%
  {\thesection.}%
  {0,3em}%
  {}%

\newcommand{\hu}[1]{ ^h\!{#1}^{ \times}}

\title{\Large\bfseries Unitarily Graded Field Extensions}
\author{\large\bfseries Holger Brenner, Almar Kaid
and Uwe Storch\\\\}

\date{}

\begin{document}

\maketitle

\newboolean{optionalremark}

\setboolean{optionalremark}{false}

\begin{myabstract}
\noindent We introduce the universal unitarily graded $A$-algebra
$A\langle U \rangle$ for a commutative ring $A$ and an arbitrary
extension $A^\times \hookrightarrow U$ of abelian groups (where
$A^\times$ denotes the group of units of $A$) and use this concept
to give among other things simplified and concise proofs of the main
theorems of the so called co-Galois theory described in the book
``Cogalois Theory'' by T. Albu. The main tool is a generalisation of
a theorem by M. Kneser which, in our language, is a criterion for $A
\langle U \rangle$ to be a field under the obvious assumption that
the base ring $A$ is itself a field and that $A \langle U \rangle$
is separable (and algebraic) over $A$. This theorem implies also the
theorem of A. Schinzel on linearly independent roots. Since
essential group extensions $A^\times \hookrightarrow U$ play a
fundamental role we discuss examples involving the injective hull of
the multiplicative group of a field. Furthermore,we develop criteria
for Galois extensions which allow a co-Galois grading, in particular
for the cyclic case.
\end{myabstract}

\bigskip


\thispagestyle{empty}

\section{Introduction}\label{introduction}
Throughout this paper we will consider only commutative rings.
First of all we fix some notations which we will use consistently:
$\PP$ denotes the set of all prime numbers in $\NN^*= \NN
\backslash \{0\}$. For an abelian group $G$ with a
multiplicatively written operation and a prime number $p$ we
denote by $G[p^\infty]:=\{x \in G: x^{p^k}=1 \mbox{, }k \geq 1\}$
the $p$-\emph{primary component} and by $G[p]:= \{x \in G: x^p =
1\}$ the $p$-\emph{socle} of $G$. The \emph{order} of $G$ is
denoted by $|G|$ or $\ord(G)$ and its \emph{exponent} by
$\exp(G)$. The order of an element $x$ of a group is denoted by
$\ord\,x (\in \NN)$. We write $A^\times$ for the group of units of
a ring $A$ and $\mu_n(A) := \{x \in A^\times: x^n = 1\}$ for the
group of $n$-th roots of unity in $A$, $n \in \NN^*$. For a field
$K$ the group $\mu_n(K)\subseteq K^\times$ is cyclic. By $\zeta_n$
we always denote a \emph{primitive} root of unity in $K^\times$,
i.e. a root of unity of order $n$. If $K=\CC$, we denote by
$\zeta_n$ the standard root of unity $\exp (2\pi i/n)$. If $K
\subseteq L$ is an extension of fields we simply write $L|K$ and
denote by $[L:K]:= \Dim_K\,L$ the \emph{degree} of $L$ over $K$.
The Galois group $\aut_{K\mbox{-}\alg}L$ of $L|K$ is denoted by
$\gal(L|K)$.

\smallskip

In this paper $A$ denotes always a base ring,
which is not the zero ring, and $D$ denotes an abelian group with additively
written operation.

\begin{defi}
Let $B= \bigoplus_{d\in D} B_d$ be a $D$-graded $A$-algebra. Then we call $B$
\emph{unitarily $D$-graded}, if $B_0 = A$ and $B_d^\times :=
B_d \cap B^\times \neq \emptyset$ for every $d \in D$.
\end{defi}

For a unitarily $D$-graded $A$-algebra $B= \bigoplus_{d \in D} B_d$
every homogeneous component $B_d$, $d \in D$, is obviously a free
$A$-module of rank one. (Notice that in the unitarily graded case
$B_d B_e = B_{d+e}$ holds for $d,e \in D$. Hence, unitarily graded
algebras are strongly graded algebras in the sense of \cite{dade}.)
In particular, a unitarily $D$-graded $A$-algebra is a free
$A$-algebra.

\smallskip

Let $x \in B_d^\times$. Then $x^{-1} \in B_{-d}$, $B_d^\times =
A^\times x$ and $x$ is transcendental over $A$ if $d \in D$ is not a
torsion element and algebraic over $A$ with minimal polynomial
$X^{\ord d} - x^{\ord d}$ else. In particular, a unitarily graded
$A$-Algebra $B$ is integral over $A$ if and only if its grading
group is a torsion group.

\smallskip

If $D^\prime \subseteq D$ is a subgroup of $D$ then $B_{D^\prime}:=
\bigoplus_{d \in D^\prime} B_d$ is obviously a unitarily
$D^\prime$-graded $A$-subalgebra of $B$. Moreover, $B$ is unitarily
$D/D^\prime$-graded over $B_{D^\prime}$ with homogeneous components
$B_{d + D^\prime} = \sum_{d^\prime \in D^\prime} B_{d +d^\prime} =
B_d B_{D^\prime}$ and $B_{d + D^\prime}^\times = B_d^\times
B_{D^\prime}^\times$. Conversely, if $C \subseteq B$ is an
$A$-subalgebra of $B$ then one easily checks that $D_C := \{d \in D:
B_d^\times \cap C^\times \neq \emptyset\}$ is a subgroup of $D$.

\smallskip

If $B$ is unitarily $D$-graded and $D = D_1 \times D_2$  with subgroups
$D_1, D_2 \subseteq D$, then the canonical homomorphism $B_{D_1} \otimes_A B_{D_2} \rightarrow B = B_D$
is an isomorphism of $D$-graded rings. If $B$ and $B^\prime$ are unitarily $D$ and
 $D^\prime$-graded respectively then $B \otimes_A B^\prime = \bigoplus_{(d,d^\prime) \in D \times
D^\prime} B_d \otimes_A B_{d^\prime}$ is a unitary $(D\times D^\prime)$-grading of $B\otimes_A
B^\prime$.

\smallskip

Let $B$ be a unitarily $D$-graded $A$-algebra and $A \rightarrow A^\prime$ a ring homomorphism.
Then $B^\prime := B \otimes_A A^\prime$ is a unitarily $D$-graded $A^\prime$-algebra.

\begin{myexam}
The $A$-algebra $A[X]/(X^n-a)$, $a \in A^\times$, has a natural unitary $\ZZ_n$-grading. Hence,
$$A[X_1,\ldots,X_r]/(X_1^{n_1}-a_1,\ldots,X_r^{n_r} - a_r) = \bigotimes
\nolimits_{j=1}^r A[X_j]/ (X_j^{n_j}-a_j),$$ $a_1,\ldots,a_r \in
A^\times$, has a natural unitary $(\, \prod_{j=1}^r
\ZZ_{n_j})$-grading. Since any finite abelian group is a direct
sum of cyclic groups \emph{every finite unitarily graded
$A$-algebra is up to $($graded$\,)$ isomorphism of this type.}
\end{myexam}

\begin{myexam}
The group algebra $A[D] = \bigoplus_{d \in D} A\, T^d$ is obviously a
unitarily $D$-graded $A$-algebra.
\end{myexam}

We denote by ${\hu {B}}$ the homogeneous units of a graded ring $B$,
which is obviously a subgroup of $B^\times$. Two unitary gradings
are by definition \emph{essentially the same} if their groups of
homogeneous units coincide. The map $\degr: {\hu {B}} \rightarrow
D$, which maps an element $x_d \in B_d^\times$ to its degree $d$, is
a homomorphism of abelian groups. By definition of a unitarily
$D$-graded $A$-algebra we get the following:

\begin{prop}\label{impsequenz}
Let $B$ be a unitarily $D$-graded $A$-algebra. Then
$$1 \longrightarrow A^\times \longrightarrow {\hu {B}}
\stackrel{\degr} {\longrightarrow} D \longrightarrow 0$$ is an exact
sequence of abelian groups. Especially, there is a canonical
isomorphism $D \cong {\hu {B}}/A^\times$.
\end{prop}

With respect to Proposition~\ref{impsequenz}, we often identify
the groups $D$ and ${\hu {B}}/A^\times$, but continue to write the
operation in $D$ additively.

\smallskip

For an abelian group $U$ containing $A^\times$, we construct a
universal unitarily  $U/A^\times$-graded $A$-algebra in the
following way: We denote $U/A^\times$ by $D$ and write $d \in D$ for
a class $A^\times x$. We choose a system $x_d \in U$ of
representatives for the elements $d = A^\times x_d \in D =
U/A^\times$ and consider the free $A$-module
$$A\langle U \rangle := \bigoplus\nolimits_{d\in D} Ax_d$$
with $A$-basis $x_d$, $d \in D$. The product $x_d x_e$ for $d,e \in
D$ is given by the multiplication in $U$, i.e. $x_d x_e = a_{d,e}
x_{d+e}$ with $a_{d,e} \in A^\times$. It is obvious that $A \langle
U \rangle$ is a unitarily $D$-graded $A$-algebra and that $U$ can be
identified with ${\hu {A \langle U \rangle}}$ via the canonical
inclusion $\gamma:~U \rightarrow A\langle U \rangle^\times$,
$x\mapsto ax_d$, where $A^\times x= A^\times x_d$ and $x = ax_d$
with $a \in A^\times$. In particular $A\langle U \rangle_d^\times =
A^\times x_d$ and for any system $y_d \in U$, $d \in D$, of
representatives for $U/A^\times$ the elements $\gamma(y_d)$, $d \in
D$, form an $A$-basis of $A\langle U \rangle$.

\smallskip

The pair $(A\langle U \rangle,\gamma)$ has the following universal
property (which, by the way, proves its uniqueness):

\begin{prop}
Let $B$ be a $($not necessarily graded$\,)$ $A$-algebra together
with a group homomorphism $\psi: U \rightarrow B^\times$ that
coincides on $A^\times$ with the structure homomorphism of $B$. Then
there is a uniquely determined $A$-algebra homomorphism $\bar{\psi}:
A\langle U \rangle \rightarrow B$ such that $\psi = \bar{\psi} \circ
\gamma$.
\end{prop}

\begin{proof}
Because the elements $x_d$ form an $A$-basis of
$A\langle U \rangle$ we can extend the group homomorphism $\psi$ to an
$A$-module homomorphism $\bar{\psi}:~A\langle U \rangle\rightarrow B$ by
$\bar{\psi}(x_d) := \psi(x_d)$. Due
to the assumption that $\psi$ coincides on $A^\times$ with the
structure homomorphism of $B$ one easily checks that $\bar{\psi}$
is even  an $A$-algebra homomorphism.
\end{proof}

\begin{myrem}
One can define the algebra $A\langle U \rangle$ alternatively as
$A \otimes_{B[A^\times]}B[U]$, where $B \rightarrow A$ is any ring
homomorphism (and $B[U]$, $B[A^\times]$ are the group algebras).
In particular, one can set $A\langle U \rangle := A
\otimes_{\ZZ[A^\times]}\ZZ[U]$. We thank the referee for this
useful comment.
\end{myrem}

\begin{myrem}
We can interpret every unitarily graded $A$-algebra $B$ as such a
universal algebra $A \langle U \rangle$ with $U:={\hu {B}}$. So the
algebra structure of $B$ is already determined by the group
extension $A^\times \hookrightarrow {\hu {B}}$.
\end{myrem}

\begin{myrem}
It is well known that the group $\ext(D,A^\times)=\ext^1_\ZZ(D,A^\times)$
describes the isomorphy classes of exact sequences
$$1 \longrightarrow A^\times \longrightarrow U \longrightarrow D
\longrightarrow 0$$ of abelian groups. So the group
$\ext(D,A^\times)$ also classifies the isomorphy types of unitarily
$D$-graded $A$-algebras. The trivial element of $\ext(D,A^\times)$
is the direct product $A^\times \times D$ which corresponds to the
group algebra $A[D] = A \langle A^\times \times D \rangle$.
\end{myrem}

\section{Unitarily graded field extensions}\label{sec2}

The aim of this section is to give an answer to the following
natural question: For which extensions $A^\times \hookrightarrow
U$ of abelian groups is the universal algebra $A\langle U \rangle$
a field? If this is the case, necessarily $A$ itself is a field.
Therefore, we assume in this section that the base ring $A$ is a
field $K$. Furthermore we use throughout our standard notations:
For an extension $K^\times \hookrightarrow U$ of abelian groups $K
\langle U \rangle$ is the universal algebra constructed in section
\ref{introduction}. It is unitarily graded, its group ${\hu {K
\langle U \rangle}}$ of homogenous units can be identified with
$U$ and the grading group is $D:=U/K^\times$. For every unitarily
graded $K$-algebra $B$ the canonical homomorphism $K\langle {\hu
{B}} \rangle \rightarrow B$ is an isomorphism. We want to clarify
that a unitarily graded field extension $L|K$ is a \emph{Kneser
extension} as introduced in \cite[Definition 2.1.9 and Definition
11.1.1]{albu} and vice versa. Important examples of unitarily
graded field extensions are the Kummer extensions.

\begin{myexam}\label{kummerbeispiel}
We recall that a (not necessarily finite) algebraic field extension
$L|K$ is a \emph{Kummer extension}, if $L|K$ is a Galois extension
with abelian Galois group $\gal(L|K)$ and if for every finite
intermediate field $K \subseteq E \subseteq L$ the base field $K$
contains a root of unity of order $\exp(\gal(E|K))$. The last
property is fulfilled if and only if the group of all continuous
characters $\check \gal(L|K):= \Hom(\gal(L|K),\QQ/\ZZ)$ can be
identified with the group of the (continuous) characters $\gal(L|K)
\rightarrow K^\times$ with values in $K^\times$.

\begin{prop}\label{kummergraduierung}
$(1)$ Let $L|K$ be a Kummer extension with Galois group $G:=
\gal(L|K)$. For a $($continuous$)$ character $\chi: G \rightarrow
K^\times$ let $L_\chi$ denote its eigenspace $L_\chi := \{ x \in L
: \sigma(x) = \chi(\sigma) x \mbox{ \rm  for all } \sigma \in
G\}$. Then $L = \bigoplus_{\chi \in \check G} L_\chi$ is a unitary
$\check G$-grading of $L$ over $K$, $\check G = \Hom(G,K^\times)$.

\smallskip

\noindent $(2)$ Conversely, let $L = \bigoplus_{d\in D} L_d$ be a unitarily
$D$-graded field extension of $K = L_0$ and suppose that $K$ contains a
root of unity of order $n_0$ whenever $D$ contains an element of order
$n_0$. Then $L$ is a Kummer extension of $K$ with Galois group
$\check D = \Hom(D,K^\times)$, where a character
$\delta: D \rightarrow K^\times$ operates as $\delta(\sum_{d \in D} x_d)
 = \sum_{d \in D} \delta(d) x_d$. $($Here a character $\delta \in \check D$
is an \emph{arbitrary} group homomorphism $D \rightarrow K^\times$,
and the topology of $\check D$ as a profinite group is given by the
finite subgroups $D_0 \subseteq D$ with the surjections $\check
D\rightarrow \check D_0$, $\check D = \varprojlim \check D_0$.$)$
Especially, $L_d$ is necessarily the eigenspace for the character
$\chi_d: \check D \rightarrow K^\times$, $\delta \mapsto \delta(d)$,
and the given grading of $L$ can be identified with the grading of
part $(1)$. Furthermore, the only intermediate fields of $L|K$ are
the graded fields $L_{D^\prime}$, $D^\prime$ subgroup of $D$.
\end{prop}

\begin{proof}
One reduces easily both assertions to the case of a finite extension
$L|K$. For part (2) note that the grading group $D$ is necessarily a
torsion group by Proposition \ref{fieldcase} below.

\smallskip

(1) Then, by the assumption on the roots of unity in $K$, every
$K$-linear operator $\sigma \in G$ of $L$ is diagonalisable over
$K$. Since $G$ is commutative the elements of $G$ are simultaniously
diagonalisable, i.e. $L=\bigoplus_{i \in I} L_i$ with $G$-invariant
$1$-dimensional $K$-subspaces $L_i \subseteq L$. Trivially, for
every $i \in I$ the function $\chi:\, G \rightarrow K^\times$ with
$\chi(\sigma)=\sigma(x) x^{-1}$ for all $\sigma \in G$ and all $x\in
L_i \backslash \{0\}$ is a character. Because of $|\check G| = |G| =
[L:K]$, and $L_\chi L_{\chi^\prime} \subseteq L_{\chi \chi^\prime}$,
it suffices to show that $\Dim_K L_\chi \leq 1$ for all $\chi \in
\check G$, but $L_1 = K$ for the trivial character $1$ and $L_\chi =
L_1x$ for any $x \in L_\chi \backslash \{0\}$.

\smallskip

(2) Obviously, $\delta:\,L\rightarrow L$ is a
$K$-automorphism of $L$, which respects the grading. Because of $|D|
= |\check D| = [L:K]$ these are all $K$-automorphisms of $L$.
\end{proof}

Let us mention that a Kummer extension $L|K$ may have unitary
gradings which are essentially different from the canonical grading
described in Proposition \ref{kummergraduierung}. For instance, the
cyclotomic field $\QQ[\zeta_8] = \QQ[i,\sqrt{2}] \cong \QQ[X]/(X^4 +
1) \cong \QQ[Y,Z]/(Y^2 +1,Z^2 -2)$ is a Kummer extension of $\QQ$
which has besides the canonical $\ZZ_2 \times \ZZ_2$-grading a
unitary $\ZZ_4$-grading. The canonical grading of a Kummer extension
$L|K$ is characterised by the property that the base field $K$
contains a root of unity of order $n_0$ if the grading group $D$
contains an element of order $n_0$, $n_0 \in \NN^*$.
\end{myexam}

\begin{prop}\label{fieldcase}
Let $L = K \langle U \rangle$ be a field. Then the group extension
$K^\times \hookrightarrow U$ is essential and, in particular, the
grading group $D=U/K^\times$ is a torsion group.
\end{prop}

\begin{proof}
To prove that $D$ is a torsion group let $d_0 \in D$, $d_0\neq 0$,
and $x_{d_0} \in L_{d_0}^\times$. Then $1+x_{d_0} \in L^\times$. Let
$\sum_{d \in D} y_d$ be the inverse of $1+x_{d_0}$. The equation
$(1+x_{d_0}) \sum_{d\in D}y_d = 1$ implies $y_0 = 1 -
x_{d_0}y_{-d_0}$ and $y_d = -x_{d_0}y_{d-d_0}$ for all $d \neq 0$.
The first equation implies $y_0 \neq 0$ or $y_{-d_0} \neq 0$. The
other equations imply (by induction) $y_{kd_0} = (-1)^k x_{d_0}^k
y_0$ for all $k \in \ZZ$, hence $y_{kd_0} \neq 0$ for all $k \in
\ZZ$. It follows that $\ZZ d_0$ is a finite group.

\smallskip

We want to recall that an extension $H \subseteq G$ of abelian
groups is by definition \emph{essential}, if for every subgroup $F
\subseteq G$ with $F \cap H = 1$ already $F = 1$ holds. It is easy
to prove that this is equivalent to the following conditions: The
quotient $G/H$ is a torsion group and, for every prime number $p$,
the $p$-socles $H[p]$ and $G[p]$ coincide. In our case $H =
K^\times$ is the multiplicative group of the field $K$. Therefore,
the extension $K^\times \subseteq U$ is essential if and only if
$U/K^\times$ is a torsion group and every root of unity of order
$p$, $p \in \PP$, in $U$ belongs already to $K^\times$.

\smallskip

The quotient $U/K^\times = D$ is a torsion group by the first part.
Assume $\zeta_p$ is a root of unity of order $p$, $p \in \PP$, in
${\hu {L}} \backslash K^\times$. Then the graded $K$-subalgebra
$K[\zeta_p] \cong K[X]/(X^p - 1)$ is not a field. Contradiction.
\end{proof}

Proposition \ref{fieldcase} says in particular, that a unitarily
graded field extension $L|K$ is algebraic. A homogeneous element
$x_d \in L_d^\times$, $d \in D \cong {\hu {L}/K^\times}$, has
degree $\ord d$ over $K$. Therefore, $L$ is separably algebraic if
and only if $\Char K = 0$ or $\Char K = \ell >0$ and
$D[\ell^\infty] = 0$.

\smallskip

Since we are only interested in the separable case, from now on
\emph{we presuppose in this section that $U/K^\times$ is a torsion
group and that $(U/K^\times)[\ell^\infty] = 1$ in case $\Char K =
\ell> 0$.}

\smallskip

The following three lemmas are the essential steps for the proof
of the main theorem.

\begin{lem}\label{pcase}
Let $D= U/K^\times$ be a finite $p$-group of order $p^\alpha$,
$\alpha \geq 1$, $p$ prime $(\neq \Char K)$. In case $p = 2$
assume $i= \sqrt{-1} \in K$. Then $B:= K \langle U \rangle$ is a
field if and only if the group extension $K^\times \hookrightarrow
U$ is essential. -- In this case $(B^\times/K^\times)[p^\infty] =
U/K^\times = {\hu {B}}/K^\times = D$.
\end{lem}

\begin{proof}
By Proposition \ref{fieldcase} the extension $K^\times
\hookrightarrow U$ is essential if $B$ is a field. For the proof of
the converse and the supplement we use induction on $\alpha$. Let
$\alpha = 1$. Then $B = K[x] \cong K[X]/(X^p-a)$ where $x \in
U\backslash K^\times$ and $a = x^p \in K^\times$. We have to show
that the polynomial $X^p - a$ is irreducible. Assume that $X^p -a$
has a zero $y$ in a field extension $L$ of $K$ of degree $m < p$.
Then $a=y^p$ and $a^m =\norm^L_K(a) = \norm_K^L(y)^p$ (where
$\norm_K^L$ denotes the norm function). Because of $\gcd(m,p) =1$ we
have $a = b^p$ with $b \in K^\times$ and $(x/b)^p =1$ with $x/b \in
U$. It follows $x/b \in K^\times$ (since $K^\times \hookrightarrow
U$ is essential) and $x \in K^\times$. Contradiction.

\smallskip

To prove the supplement it is enough to show: If $y \in B^\times$
and $y^p \in U= {\hu {B}}$ then $y \in U$. We adjoin if necessary to
$K$ a root of unity $\zeta_p$ of order $p$ and consider the Kummer
extension $K[\zeta_p] \subseteq K[\zeta_p] \otimes_K B = B[\zeta_p]
\cong K[\zeta_p][X]/(X^p - a)$. (Note that $K[\zeta_p] \otimes B$ is
a field because of $\gcd([K[\zeta_p]:K],[B:K]) =1$.)

\smallskip

First assume that even $y^p \in K^\times$. If $y \notin K^\times$
then $B=K[y]$ and $B[\zeta_p]= K[\zeta_p][y]$. By Proposition
\ref{kummergraduierung} the element $y$ is homogeneous in
$B[\zeta_p]$ (since $K[\zeta_p]y^k$, $k=0,\ldots,p-1$, are the
homogeneous components of a unitary grading of $B[\zeta_p]$). Then
$y$ is homogeneous in $B$ too, i.e. $y \in U$.

\smallskip

Now suppose $y^p \notin K^\times$. Then $y^{p^2}=(y^p)^p =:c \in K^\times$
and $X^p-c$ is the minimal (= characteristic) polynomial of $y^p$ and
$c=(-1)^{p+1} \norm_K^B (y^p) = (-1)^{p+1} \norm_K^B(y)^p$. In any case
$c$ is a $p$-th power in $K^\times$ (in case $p =2$ we use $i \in K$).
This contradicts the irreducibility of $X^p - c$.

\smallskip

For the induction step assume $|D| = p^{\alpha +1}$. Let
$\tilde{D} \subset D$ be a subgroup of order $p^\alpha$. Then by
induction hypothesis, the unitarily $\tilde{D}$-graded subalgebra
$\tilde{B}:=B_{\tilde{D}} \subset B$ is a field with
$(\tilde{B}^\times/K^\times)[p^\infty]={\hu {\tilde{B}
}}/K^\times$ and $B$ is a unitarily $D/\tilde{D}$-graded
$\tilde{B}$-algebra with $\tilde{B}^\times {\hu {B}}$ as group of
homogeneous units. The group extension $\tilde{B}^\times
\hookrightarrow \tilde{B}^\times{\hu {B}}$ is essential. To prove
this, let $(yz)^p = y^pz^p =1$, $y \in \tilde{B}^\times$, $z \in
{\hu {B}}$. Then $z^p \in {\hu {\tilde{B}}}$, $y^p \in {\hu
{\tilde{B}}}$, so $y \in {\hu {\tilde{B}}}$ by the induction
hypothesis on the supplement. Hence $yz \in {\hu {B}}$ and $yz \in
K^\times \subseteq \tilde{B}^\times$ since $K^\times
\hookrightarrow {\hu {B}}$ is essential. The case $\alpha =1$
implies that $B$ is a field and
$(B^\times/\tilde{B}^\times)[p^\infty] = \tilde{B}^\times {\hu
{B}}/\tilde{B}^\times$.

\smallskip

To prove $(B^\times/K^\times)[p^\infty] = {\hu {B}}/K^\times$ let
$w \in B^\times$ represent an element in
$(B^\times/K^\times)[p^\infty]$. Then $w \in \tilde{B}^\times {\hu
{B}}$, $w = uv$ with $u \in \tilde{B}^\times$, $v \in {\hu {B}}$,
hence $u \in {\hu {\tilde{B}}}$ and $w \in {\hu {B}}$ as wanted.
\end{proof}

\begin{lem}\label{2caseone}
Let $D = U/K^\times$ be a finite $2$-group of order $2^\alpha$,
$\alpha \geq 1$. Assume $U$ contains no element of order $4$. Then
$B:=K \langle  U \rangle$ is a field if and only if the group
extension $K^\times \hookrightarrow U$ is essential. -- In this
case $(B^\times/K^\times)[2^\infty] = U/K^\times = {\hu
{B}}/K^\times = D$.
\end{lem}

\begin{proof}
By Proposition \ref{fieldcase} the extension $K^\times
\hookrightarrow U$ is essential if $B$ is a field. We consider the
extension $K[i] \subseteq B[i]:= K[i] \otimes_K B$. It is enough
to show that the extension $K[i]^\times \hookrightarrow {\hu
{B[i]}}$ is essential. Then, due to \ref{pcase}, $B[i]$ is a field
hence $B$ too. Furthermore,
$(B[i]^\times/K[i]^\times)[2^\infty]={\hu {B[i]}} /K[i]^\times$
which implies $(B^\times/K^\times)[2^\infty] = {\hu {B}/K^\times}$
because of ${\hu {B}}= {\hu {B[i]}} \cap B$. We have $B[i]_d = B_d
\oplus B_d i$ for all $d \in D$. So let $b,c \in B_d$ with $1=
(b+ci)^2 = b^2 + 2bci - c^2$. Comparison of coefficients yields
$b^2 - c^2 = 1$ and $2bc = 0$. Because $\Char K \neq 2$ we have $b
= 0$ or $c = 0$. Suppose $b = 0$, hence $-c^2 = 1$. But this means
$c = \pm i \in {\hu {B}}$ which is a contradiction. So we have $c
= 0$, hence $b^2 = 1$. Because $K^\times \subseteq {\hu {B}}$ is
essential we get $b = \pm 1$.
\end{proof}

Note that in the situation of Lemma \ref{pcase} or Lemma
\ref{2caseone} the torsion group $\tor(B^\times/K^\times)$ may be
larger than $U/K^\times = {\hu {B}}/K^\times$ even if $B$ is a
field! A simple example is $B = \QQ[\zeta_3]= \QQ[\sqrt{-3}]$ over
$K = \QQ$.

\smallskip

If $D = U/K^\times$ is a finite $2$-group the condition that the
extension $K^\times \hookrightarrow U$ is essential is in general
not sufficient for $K \langle U \rangle$ to be a field. By
\ref{2caseone} this can only occur if $U$ contains an element of
order $4$.

\vskip 0,2 cm

\begin{myexam}\label{2examone}
We consider the polynomial $X^4 + 4 \in \QQ[X]$. We have the
well-known decomposition $X^4 + 4 = (X^2 - 2X +2) ( X^2 +2X +2)$
over $\QQ$, so the unitarily $\ZZ_4$-graded $\QQ\,$-Algebra $B :=
\QQ[X]/(X^4 + 4)$ is not a field. But the extension $\QQ^\times
\subseteq {\hu {B}}$ is essential due to the fact that there is no
element $y \in {\hu {B}}\backslash \QQ^\times$ with $y^2 = 1$. The
element $x^2/2$ has order $4$ in $U = {\hu {B}}$.
\end{myexam}

\vskip 0,1 cm

\begin{lem}\label{2casetwo}
Let $D =U/K^\times$ be a finite $2$-group of order $2^\alpha$,
$\alpha \geq 1$. Assume $U$ contains an element of order $4$ which
is not an element of $K$. Then $B := K \langle U \rangle$ is a
field if and only if the group extension $K^\times \hookrightarrow
U$ is essential and $-4 \notin U^4$ $($i.e. there is no element
$x\in U$ with $x^4 = -4)$.
\end{lem}

\begin{proof}
If $B$ is a field then $K^\times \hookrightarrow U$ is essential by
Proposition \ref{fieldcase}. Furthermore, if there is an element $x \in U$
with $x^4 = -4$, then $x$ represents an element of order $4$ in $D =
U/K^\times$ because $(x^2/2)^2 = -1$ and therefore $x^2 = \pm 2i
\notin K^\times$ by assumption. It follows $K[x] \cong K[X]/(X^4 + 4)$, and
$K[x]$ is not a field because of $X^4 + 4 = (X^2-2X +2)(X^2+2X +2)$, see also
Example \ref{2examone}.

\smallskip

Conversely, the element $i \in U$ of order $4$ represents an element
of order $2$ in $U/K^\times$ because of $i^2 = -1$. So $K[i]
\subseteq B$ is a graded quadratic subfield of $B$ and $B$ is
unitarily graded over $K[i]$ with grading group $K[i]^\times
U/K[i]^\times$. By Lemma \ref{pcase} it is now sufficient to show
that the extension $K[i]^\times \hookrightarrow K[i]^\times U$ is
essential. To do this, let $y^2 = x^2$ with $y\in U$, $x = a+bi \in
K[i]^\times$, $a,b \in K$. Then $y^2 = a^2 - b^2 + 2abi \in U$,
hence $a^2-b^2 = 0$ or $2ab =0$. If $a^2 - b^2 = 0$, then $a = \pm
b$, $(y/a)^4 =(\pm2i)^2 = -4$ which is impossible by assumption.
Therefore $ab = 0$, i.e. $x \in U$, hence $x^{-1}y \in U$ and
$x^{-1}y = \pm 1$ since $K^\times \hookrightarrow U$ is essential.
\end{proof}

\begin{myrem}
(1) In the situation of \ref{2casetwo} it is rather difficult to
describe the $2$-torsion group $(B^\times/K^\times)[2^\infty]$.
Because $1+i \notin U$ represents an element of order $4$ in
$B^\times/K^\times$ the group $(B^\times/K^\times)[2^\infty]$ is
always larger than ${\hu {B}}/K^\times =U/K^\times$. But the
simple example $K:= \RR$, $B:= \RR[i] = \CC$ shows that
$(B^\times/K^\times)[2^\infty]$ can be much larger than
$U/K^\times$.

\smallskip

\noindent (2) It would be interesting to understand the structure
of the separable $K$-algebra $B = K \langle U \rangle$ or at least
its spectrum if the essential extension $K^\times \hookrightarrow
U$ fulfills all the assumptions of Lemma \ref{2casetwo} and
moreover $-4 \in U^4$. For illustrations look at Example
\ref{2examone} and its extension Example \ref{injhullqexam} in the
next section or at the following one: For $K$ take the real number
field $\QQ[\zeta_{16}]\cap \RR$ and for $U$ the essential
extension $K^\times \mu_{16}(\CC)$ of $K^\times$ with $K^\times
\mu_{16}(\CC)/K^\times \cong \ZZ_8$. Then $1+i =\sqrt{2}\,\zeta_8
\in U$ with $(1+i)^4 = -4$ and $K\langle U \rangle \cong K
\otimes_\QQ \QQ[\zeta_{16}]$ splits into $4$ components which are
isomorphic quadratic field extensions of $K$.

\smallskip

The comments in this remark also show that the statements in
\cite{schejastorch2}, \S 93, Exercise 14e)(3),(4) are not correct.
\end{myrem}

The following theorem which generalises amongst others the Theorem
of M. Kneser in \cite{kneserlinearwurzel} is the main result and
summarises the results of the previous lemmata, cf. also
\cite[Satz 3.2.6]{kaiddiplom}.

\begin{theo}\label{maintheo}
For the group extension $K^\times \hookrightarrow U$ $($with
$(U/K^\times)[\ell^\infty] = 1$ if $\Char K = \ell > 0)$ the
universal algebra $K\langle U \rangle$ is a field if and only if
the extension $K^\times \hookrightarrow U$ is essential and
moreover $-4 \notin U^4$ in case $U$ contains an element of order
$4$ not in $K^\times$. -- In this case $(K\langle U
\rangle^\times/K^\times)[p^\infty] = U/K^\times$ if $U/K^\times$
is a $p$-group, $p \geq 3$, and $(K\langle U
\rangle^\times/K^\times)[2^\infty] = U/K^\times$ if $U/K^\times$
is a $2$-group and $U$ contains no element of order $4$ not in
$K$.
\end{theo}

\begin{proof}
Let $D:= U/K^\times$. If the unitarily $D$-graded $K$-algebra $B
:= K\langle U \rangle$ is a field then $K^\times \hookrightarrow
U$ is essential by \ref{fieldcase} and the exceptional case is
settled by Lemma \ref{2casetwo} because of $B_{D[2^\infty]}
\subseteq B$.

\smallskip

Conversely, let $K^\times \hookrightarrow U$ be essential with $-4
\notin U^4$ in the special case. Because of $K\langle U \rangle =
\varinjlim K\langle U^\prime \rangle$ where $U^\prime$ runs
through the subgroups $U^\prime \subseteq U$ with $K^\times
\subseteq U^\prime$ and finite index $[U^\prime : K^\times]$ we
may assume that $D = U/K^\times$ is finite. Then $B =
\bigotimes_{p} B_{D[p^\infty]}$ because of $D = \bigoplus_{p}
D[p^\infty]$, where $p$ runs through the prime divisors of $|D|$.
Since the dimensions $\Dim_K B_{D[p^\infty]}$ are pairwise coprime
it is enough to show that all the $K$-algebras $B_{D[p^\infty]}$
are fields. But $B_{D[p^\infty]} = K \langle {\hu {
B_{D[p^\infty]}}} \rangle$ and ${\hu { B_{D[p^\infty]}}} \subseteq
{\hu {B}} = U$ are essential extensions of $K^\times$ such that
$[{\hu { B_{D[p^\infty]}}} : K^\times]$ is a power of $p$ and the
results follow from Lemma \ref{pcase}, Lemma \ref{2caseone} and
Lemma \ref{2casetwo}.
\end{proof}

If the factor group $U/K^\times$ of the extension $K^\times \hookrightarrow U$
is a finite cyclic group Theorem \ref{maintheo} is the well known Theorem of Capelli
(for the separable case).

\smallskip

Obviously, if $K \langle U \rangle$ is a field then $K \langle U
\rangle$ is a Galois extension of $K$ if and only if the grading
group $D = U/K^\times$ has the following property: if $D$ contains
an element of order $n_0$ then $K \langle U \rangle$ contains a root
of unity of order $n_0$. (Note that $K \langle U \rangle$ is by our
general assumption always separable.)

\section{Applications and Examples}

In this section we prove some consequences of the results of
section \ref{sec2}. First of all we mention the following slight
generalisation of the theorems of Kneser and Schinzel in
\cite{kneserlinearwurzel} and \cite[Theorem 1]{schinzel}, see also
\cite[Theorem 2.2.1 and Theorem 11.1.5]{albu},\cite[Theorem
1.12]{stefan} and \cite[\S 93, Exercise 14]{schejastorch2}.

\begin{theo}\label{kneschi}
Let $L|K$ be a field extension with
$(L^\times/K^\times)[\ell^\infty] = 1$, i.e. ${L^\times}^\ell \cap
K^\times = {K^\times}^\ell$, if $\Char K = \ell > 0$, and let $U
\supseteq K^\times$ be a subgroup of $L^\times$. Furthermore, let
$x_i$, $i \in I$, be a full system of representatives for the
elements of $U/K^\times$. Then $E:= \sum_{i \in I} Kx_i$ is a
$K$-subalgebra of $L$ and the following conditions are equivalent:
\begin{enumerate}
\item $E$ is a field and the $x_i$, $i \in I$, are linearly independent over $K$.
\item $K^\times \hookrightarrow U$ is an essential extension of groups and $1+i \notin U$
if $U$ contains a root of unity $i$ of order $4$ not in $K^\times$.
\end{enumerate}
If these conditions hold $E$ is a separable algebraic field extension of degree $[E:K]
= [U:K^\times]$.
\end{theo}

\begin{proof}
First of all, the extension $K^\times \subseteq U$ fulfills by
assumption the condition $(U/K^\times)[\ell^\infty] = 1$ if $\Char
K = \ell> 0$. Consider the universal algebra $K \langle U \rangle$
and the canonical $K$-algebra homomorphism $\psi:\, K \langle U
\rangle \rightarrow E$ induced by the inclusion $U \rightarrow
E^\times$. Condition (1) is equivalent with the condition that $K
\langle U \rangle$ is a field. Now apply Theorem \ref{maintheo}.
\end{proof}

Note that in \ref{kneschi} the algebra $E$ is a priori a field if the extension $L|K$
is algebraic.

\smallskip

The following definitions and results are inspired by the book
\cite{albu} of T. Albu and the article \cite{greitherharrison} of C.
Greither and D. K. Harrison. We also mention the work \cite{stefan}
of D. Stefan where one can find similar graded formulations for
finite field extensions.

\begin{defi}
A group extension $K^\times \hookrightarrow U$ with factor group $D = U/K^\times$ and
universal unitarily $D$-graded $K$-algebra $L:= K \langle U \rangle$ is called
\emph{co-Galois} if the following conditions are fulfilled:
\begin{enumerate}
\item $L$ is a field and $D[\ell^\infty] =0$ if $\Char K = \ell >
0$. \item Every intermediate field $K \subseteq E \subseteq L$ is
graded, i.e. $E= L_{D^\prime}$ for some subgroup $D^\prime
\subseteq D$.
\end{enumerate}
\end{defi}

\vskip 0,2 cm

We call a field extension $L|K$ \emph{co-Galois} if there exists a
co-Galois group extension $K^\times \hookrightarrow U$ such that
$L \cong K \langle U \rangle$. In this case the extension
$K^\times \subseteq U$ is uniquely determined as we will see after
the proof of Theorem \ref{cogaloiscrit}, therefore we drop $U$
from our notation. The condition $D[\ell^\infty] = 0$ if $\Char K
= \ell > 0$ implies that a co-Galois extension is a separable
(algebraic) field extension. A co-Galois extension $L|K$ is our
graded equivalent of a \emph{U-Cogalois} extension introduced in
\cite[Definition 4.3.3 and Definition 12.1.1]{albu}.

\smallskip

For a co-Galois extension $K \subseteq L=K \langle U \rangle$ and
a subgroup $D^\prime \subseteq D = U/K^\times$ the subfield
$L_{D^\prime}$ is co-Galois over $K$ and $L$ is co-Galois over
$L_{D^\prime}$ (with respect to the induced $D/D^\prime$-grading).
We have maps $D^\prime \mapsto L_{D^\prime}$ and $E \mapsto D_E$
between the set of subgroups of $D$ and the set of intermediate
fields of $L|K$ which are inverse to each other. Hence, they are
(lattice) isomorphisms.

\smallskip

If $L = K \langle U \rangle$ is co-Galois and $x = \sum_{d\in D}
x_d$ is an element in $L$ then $K[x] = K[x_d: d \in D] =
L_{\langle \supp x \rangle}$ where $\langle \supp x \rangle$ is
the subgroup of $D$ generated by the \emph{support} $\supp x :=
\{d \in D: x_d \neq 0\}$ of $x$. In particular, $[K[x]:K] =
|\langle \supp x \rangle|$ and $K[x] = L$ if and only if $\langle
\supp x \rangle = D$ (cf. also \cite[Theorem 8.1.2 and Proposition
10.1.12]{albu} and \cite[Proposition 2.6]{stefan}). \emph{If $L$
is co-Galois then any $x \in L^\times$ with $x^2 \in K^\times$ is
homogeneous. Proof:} If $x \notin K$ then $[K[x]:K] = 2$, $\Char K
\neq 2$ and $x = x_0 + x_d$ with $2d = 0$, $x^2 = x_0^2+x_d^2+2x_0
x_d= x_0^2 + x_d^2$ implies $x_0 x_d =0$, i.e. $x_0 = 0$. Examples
of co-Galois extensions are the Kummer extensions, cf. Proposition
\ref{kummergraduierung}.

\smallskip

For the following characterisation of co-Galois extensions compare
also \cite[Theorem 4.3.2]{albu} and  \cite[Theorem 2.5]{stefan}
for the case of an finite extension and \cite[Theorem
12.1.4]{albu} for the infinite case.

\begin{theo}\label{cogaloiscrit}
The group extension $K^\times \hookrightarrow U$ with factor group
$D= U/K^\times$ and universal unitarily $D$-graded $K$-algebra $L:=
K\langle U \rangle$ is co-Galois if and only if the following conditions
are fulfilled:
\begin{enumerate}
\item $D$ is a torsion group with $D[\ell^\infty] = 0$ if $\Char K
= \ell> 0$. \item For all primes $p$ with $D[p^\infty] \neq 0$
every element of order $p$ in $L^\times$ belongs to $K^\times$.
\item If $D$ and $K\langle U \rangle^\times$ contain elements of
order $4$ then $K^\times$ contains an element of order $4$.
\end{enumerate}
\end{theo}

\begin{proof}
Let $L = K\langle U \rangle$ be co-Galois. Then $K^\times \hookrightarrow U$
is essential by \ref{fieldcase} and, in particular, $D = U/K^\times$ is
a torsion group.

\smallskip

Assume now that $D$ contains an element of prime order $p$ and let
$x \in U$ represent such an element. Furthermore, let $\zeta_p \neq 1$
a $p$-th root of unity in $L$. Then $\prod_{k=0}^{p-1}(X- \zeta_p^k x)
= X^p - x^p$ is the minimal polynomial over $K$ for all the elements
$\zeta_p^k x$, $k = 0,\ldots, p-1$. The subfield $K[x,\zeta_p]$ is of degree
$pm$ with $m < p$ and hence contains only one subfield of degree $p$ over $K$
since all subfields are graded. It follows $K[x] = K[\zeta_px]$ and
$\zeta_p = (\zeta_px)/x \in K[x]$, i.e. $\zeta_p \in K$.

\smallskip

Let $i \in L^\times$ be a root of unity of order $4$ and let $x \in U$ be an
element representing an element of order $4$ in $D$. Then $i$ is
homogeneous and $\prod_{k=0}^3
(X- i^k x) = X^4 - x^4$ is the minimal polynomial over $K$ for all the
elements $i^kx$, $k=0,\ldots,3$. Furthermore, $((1+i)x)^4 = (x+ ix)^4 =-4x^4
\in K^\times$, hence $[K[(1+i)x]:K] \leq 4$. If $i \notin K^\times$ then
$ix$ is homogeneous with $\deg x \neq \deg ix$ and therefore
$K[(1+i)x]=K[x,ix] = K[x] = K[ix]$ and $i \in K[x]$, $K[i] = K[x^2]$, i.e.
$\deg i = \deg x^2 = 2 \deg x$, which implies $((1+i)x)^2 = 2ix^2 \in K^
\times$. This is a contradiction!

\smallskip

To prove that conversely conditions (1),(2),(3) imply that $L =
K\langle U \rangle$ is co-Galois over $K$ we can assume that $D =
U/K^\times$ is finite.

\smallskip

Conditions (1) and (2) imply that the extension $K^\times
\hookrightarrow U$ is essential. Suppose that $U$ contains an
element $y$ of order $4$ not in $K^\times$, and assume that $x^4 =
-4$, $x \in U$. This implies $y^2 = -1 =(x^2/2)^2$, hence $y = \pm
x^2/2$ (since $K^\times \hookrightarrow U$ is essential) and $x^2
\notin K^\times$. Therefore, $x$ represents an element of order $4$
in $D$. By assumption (3) this implies that $K^\times$ contains an
element $i$ of order $4$. Then $(y/i)^2=1$ and $y/i = \pm 1$, $y=
\pm i \in K^\times$. Contradiction. By Theorem \ref{maintheo} $L$ is
a field.

\smallskip

Now, let $E$ be an intermediate field, $K \subseteq E \subseteq L =
K\langle U \rangle$. We have to show: $E= K \langle U\cap E^\times
\rangle$. Consider the group extension $E^\times \hookrightarrow
E^\times U \, \,(\subseteq L^\times)$ with index $[E^\times U :
E^\times] = [U : U\cap E^\times ]$. If the universal algebra
$E\langle E^\times U \rangle$ is a field then the canonical
homomorphism $E \langle E^\times U \rangle \rightarrow L =
E[E^\times U]$ is an isomorphism which implies $[L:E] = [E^\times U:
E^\times] = [U : U \cap E^ \times] = [L: K\langle U \cap E^\times
\rangle]$ and $E= K \langle U \cap E^\times \rangle$ because of $K
\langle U \cap E^\times \rangle \subseteq E$.

\smallskip

So we have to verify that $E^\times \hookrightarrow E^\times U$
fulfills the conditions of Theorem \ref{maintheo}. The assumption
(2) implies that $E^\times\hookrightarrow E^\times U$ is essential.
Now suppose that $E^\times U$ contains an element $i$ of order $4$
not in $E^\times$ and $x^4 = -4$ with $x \in E^\times U$. The
element $x$ represents an element of order $4$ in $E^\times
U/E^\times \cong U/U\cap E^\times$ because $(x^2/2)^2 = -1 = i^2$
and $x^2 = \pm 2i \notin E^\times$. But then $D = U/K^\times$
contains an element of order $4$ and by condition (3) $i \in K$.
Contradiction.
\end{proof}

We remark \emph{that for a co-Galois extension $L = K \langle U
\rangle$ of $K$ the group $U = \hu {L}$ of homogeneous units is
uniquely determined,} cf. also \cite[Corollary 4.4.2 and Corollary
10.1.11]{albu}. ($L$ may have however unitary gradings which are
not co-Galois, cf. Example \ref{kummerbeispiel}.) \emph{Proof:}
Let $L = K\langle U^\prime \rangle$ be another co-Galois grading
and let $x \in U^\prime$. We have to show $x \in U$. We may assume
that the order of $x$ in $U^\prime/ K^\times$ is a power of a
prime $p$, i.e that $[K[x] : K] = p^\alpha$, $\alpha \geq 1$, and
that $L = K[x]$. If $p\geq 3$ then $x$ represents an element of
$(L^\times/K^\times)[p^\infty]$ and belongs therefore to $U$ by
Theorem \ref{maintheo}.

\smallskip

If $p= 2$ then again $x \in U$. This follows from \ref{maintheo}
if $U$ does not contain an element of order $4$ not in $K^\times$.
If $i = \sqrt{-1} \in U$, $i \notin K^\times$, then $D$ is an
elementary abelian $2$-group by condition (3) in Theorem
\ref{cogaloiscrit} and the homogeneous elements $x \in L$ for both
gradings are characterised by the condition $x^2 \in K$ (cf. also
Proposition \ref{kummergraduierung}).\hspace{5,4cm}\qedsymbol

\smallskip

Furthermore, \emph{if $L = K\langle U \rangle$ is a co-Galois
extension then $(L^\times/K^\times)[p^\infty] =
(U/K^\times)[p^\infty]$ for every prime $p \geq 3$ with
$(U/K^\times)[p^\infty] \neq 1$ and the equality
$(L^\times/K^\times)[2^\infty] = (U/K^\times)[2^\infty]$ holds in
the following cases: $(1)$ $U/K^\times$ contains an element of
order $4$, $(2)$ $i (=\sqrt{-1}) \in K^\times$, $(3)$ $i \notin
L^\times$. In any case the equality $(L^\times/ K^\times)[2] =
(U/K^\times)[2]$ holds}(compare also with \cite[Theorem 4.4.1 and
Theorem 12.1.8]{albu}). The equality
$(L^\times/K^\times)[p^\infty] = (U/K^\times)[p^\infty]$ for a
prime number $p\geq 2$ is equivalent to the property that $K
\langle \Tor_p(L^\times/K^\times) \rangle$ is a field, where
$\Tor_p(L^\times|K^\times)$ is by definition the canonical
preimage of $(L^\times/K^\times)[p^\infty]$ in $L^\times$ and this
is checked with Theorem \ref{maintheo} using the characterisation
of co-Galois extensions in Theorem \ref{cogaloiscrit}.

\smallskip

Let $\Tor(L^\times|K^\times)=\{x \in L^\times:\, x^n \in K^\times
\mbox{ for some } n \} \subseteq L^\times$ denote the canonical
preimage in $L^\times$ of the torsion subgroup
$\tor(L^\times/K^\times)$ of $L^\times/K^\times$. (In
\cite{greitherharrison} the group $\tor(L^\times/K^\times) $ is
called the \emph{co-Galois group} of $L|K$.)

\begin{defi}
A field extension $L|K$ is called \emph{absolutely co-Galois} if the canonical $K$-algebra
homomorphism $K\langle \Tor(L^\times|K^\times) \rangle \rightarrow L$ induced by the
inclusion $\Tor(L^\times|K^\times) \hookrightarrow L^\times$ is an isomorphism.
\end{defi}

\smallskip

In an equivalent, but different approach finite absolutely
co-Galois extensions were treated in \cite{greitherharrison} and
called \emph{cogalois extensions}, see also \cite[Definition
12.2.1]{albu} for the infinite case. We prefer the term absolutely
co-Galois in order to stress that the grading group is the whole
torsion group of $L^\times/K^\times$.

\smallskip

If $L|K$ is absolutely co-Galois then $L$ is unitarily
$\tor(L^\times/K^\times)$-graded and ${\hu {L}} =
\Tor(L^\times|K^\times)$. \emph{The extension is necessarily
separable.} For the \emph{proof}, let $x \in L^\times$, $x^\ell
\in K^\times$, $\ell:= \Char K >0$. Then $(1+x)^\ell \in K^\times$
which implies $x \in K$ since $1,x, 1+x$ are homogeneous. This
means $(L^\times/K^\times)[\ell^\infty] = 1$.

\smallskip

The following characterisation of absolutely co-Galois extensions
is a direct consequence of Theorem \ref{maintheo}. One compares
also \cite[Theorem 1.5]{greitherharrison} and \cite[Theorem
3.1.7]{albu} for finite extensions as well as \cite[Theorem
12.2.2]{albu} for the infinite case.

\begin{theo}\label{abscogaloiscrit}
A field extension $L|K$ is absolutely co-Galois if and only if the following conditions
are fulfilled:
\begin{enumerate}
\item The group $\Tor(L^\times|K^\times)$ generates $L$ as a
$K$-algebra, the group extension $K^\times \hookrightarrow
\Tor(L^\times|K^\times)$ is essential and
$(L^\times/K^\times)[\ell^\infty] =
\Tor_\ell(L^\times|K^\times)/K^\times = 1$, i.e. ${L^\times}^\ell
\cap K^\times = {K^\times}^\ell$, if $\Char K = \ell > 0$. \item
If $L^\times$ contains a root of unity $i$ of order $4$ then $i$
belongs already to $K^\times$.
\end{enumerate}
\end{theo}

For the following two easy corollaries compare also \cite[Theorem
1.6.(a)]{greitherharrison}, \cite[Proposition 3.2.2.(2) and
Theorem 12.2.4.(4)]{albu} and \cite[Theorem 12.2.3]{albu}
respectively.

\begin{cor}
If $L|K$ is an absolutely co-Galois extension, then for any intermediate field $E$ the
extensions $L|E$ and $E|K$ are absolutely co-Galois too.
\end{cor}

Theorem \ref{cogaloiscrit} implies:

\begin{cor}
An absolutely co-Galois extension $L|K$ is co-Galois with respect to the group extension
$K^\times \hookrightarrow \Tor(L^\times|K^\times)$ and with grading group $\Tor(L^\times|K^\times)/
K^\times = \tor(L^\times/K^\times)$.
\end{cor}

Co-Galois extensions are not necessarily absolutely co-Galois.
Look at $\QQ[\zeta_3]|\QQ$ or as an extreme case $\CC|\RR$. A
co-Galois extension $L = K \langle U \rangle$ over $K$ is
absolutely co-Galois if and only if the following conditions are
fulfilled: (1) Any root of unity $\zeta_q$ of prime order $q$ in
$L^\times$ with $(U/K^\times)[q^\infty] = 1$ belongs already to
$K^\times$. (2) If the element $i$ of order $4$ belongs to
$L^\times$ then $i \in K^\times$. (If $i \notin K^\times$ then
$K[i]$ is never absolutely co-Galois.)

\begin{myexam}
Let $K$ be a field which contains for every prime $p \neq \Char K$
a root of unity of order $p$ and a root of unity of order $4$ if
$\Char K \neq 2$. Furthermore, let $\bar{K}_{\sep}$ be the
separable algebraic closure of $K$. Then the group
$\Tor(\bar{K}_{\sep}^\times|K^\times)$ is an essential extension
of $K^\times$, indeed $\Tor(\bar{K}_{\sep}^\times|K^\times)
=\I^\prime(K^\times)$ where $\I^\prime(K^\times) \subseteq
\I(K^\times)$ is the preimage of $\prod_{p \in\PP,p\neq \Char K}
(\I(K^\times)/K^\times)[p^\infty]$ in the injective hull
$\I(K^\times)$ of the group $K^\times$. The equality
$\I^\prime(K^\times) = \I(K^\times)$ holds if and only if $K$ is a
perfect field.

\smallskip

Since the group extension $K^\times \hookrightarrow
\Tor(\bar{K}_{\sep}^\times|K^\times) = I^\prime(K^\times)$ is co-Galois by
Theorem \ref{cogaloiscrit} the canonical homomorphism
$K\langle I^\prime(K^\times) \rangle \rightarrow \bar{K}_{\sep}$ is injective
and \emph{its image $K[I^\prime(K^\times)]$ is the largest absolutely
co-Galois extension of $K$}, cf. Theorem \ref{abscogaloiscrit}. It is also a
Galois extension which contains \emph{all} roots of unity, i.e. for any $n \in \NN^*$
with $n \neq 0$ in $K$ there is a root of unity of order $n$ in $K[\I^\prime
(K^\times)]$.

\smallskip

Furthermore, if $K$ contains \emph{all} roots of unity then this
extension coincides with the largest Kummer extension of $K$ which
is in this case also the largest abelian extension $\bar{K}_{\ab}$
of $K$. The Galois group of this extension is the cha\-racter group
$\Hom(\I^\prime(K^\times)/K^\times,K^\times) =
\Hom(\I^\prime(K^\times)/K^\times,\QQ/\ZZ)$ of
$\I^\prime(K^\times)/K^\times =
\tor(\bar{K}_{\sep}^\times/K^\times$), cf. Proposition
\ref{kummergraduierung}.

\smallskip

So if we iterate this construction starting with
$K_1:=K[\I^\prime(K^\times)]$ instead of $K_0:= K$ we get the Kummer
extension $K_2:= K_1[\I^\prime(K_1^\times)]$ of $K_1$ and altogether
a tower of subfields $K = K_0 \subseteq K_1 \subseteq K_2 \subseteq
\ldots$ of $\bar{K}_{\sep}$ such that every extension $K_{j+1}|K_j$,
$j \in \NN$, is absolutely co-Galois (and Kummer for $j >0$).

\smallskip

If $F$ is an arbitrary field then take for $K_0$ the field $K:=
F[\zeta_p, p \in \PP, p \neq \Char F; i]$, where $\zeta_p \in
\bar{F}_{\sep} (= \bar{K}_{\sep}$) is a root of unity of order $p$
(and $i \in \bar{K}_{\sep}$ of order $4$ if $\Char F \neq 2$). If
$\Char F = 0$ then $\bigcup_{j \geq 0} K_j = \bar{F}_{\solv}$ is the
union of all Galois extensions of $F$ in $\bar{F}_{\sep} = \bar{F}$
with solvable Galois group.
\end{myexam}

\begin{myexam}
Let $K$ be an \emph{ordered} field and let $\bar{K}_{\real}$ be the
real closure of $K$. Then the group
$\Tor(\bar{K}_{\real}^\times|K^\times)$ is an essential extension of
$K^\times$ since $\pm 1$ are the only roots of unity in
$\bar{K}_{\real}^\times$, indeed
$\Tor(\bar{K}_{\real}^\times|K^\times) = \{\pm 1\} \I(K_+^\times)$,
where $\Tor(\bar{K}_{\real,+}^\times|K_+^\times) = \I(K_+^\times)
\subseteq \bar{K}_{\real,+}^\times$ is the injective hull of the
group of positive elements in $K$.

\smallskip

Since the group extension $K^\times \hookrightarrow \Tor(\bar{K}_{\real}^\times|K^\times)$ is
co-Galois by Theorem \ref{cogaloiscrit} \emph{the canonical homomorphism
$K\langle \{\pm 1\} \I(K_+^\times) \rangle \rightarrow \bar{K}_{\real}$ is injective} and
\emph{its image $K[\I(K_+^\times)]$ is the largest co-Galois extension of $K$ in
$\bar{K}_{\real}$. It is even absolutely co-Galois}, cf. Theorem \ref{abscogaloiscrit}.

\smallskip

In case that $K = \QQ$ or, more generally, that $K$ is a real
algebraic number field the injectivity of the canonical map
$K\langle \{\pm1\} \I(K_+^\times) \rangle \rightarrow
\bar{K}_{\real} \subseteq \RR$ is a classical result of Besicovitch
\cite{besicovitch} and Siegel \cite{siegelwurzel}.

\smallskip

That $\QQ \subseteq \QQ[\I(\QQ_+^\times)]$ is a co-Galois
extension can be expressed in the following way: If
$(\nu_{1\sigma},\ldots,\nu_{r\sigma}) \in \QQ^r$, $\sigma =
1,\ldots,s$, are $r$-tuples which represent different elements in
$(\QQ/\ZZ)^r$ and if $p_1,\ldots,p_r$ are different prime numbers
then the degree of every element
$$x = \sum\nolimits_{\sigma=1}^s a_\sigma p_1^{\nu_{1\sigma}} \cdots p_r^{\nu_{r\sigma}}$$
with $a_1, \ldots, a_s \in \QQ^\times$ is $|d|$ where $1/d$ is the
greatest common divisor of \emph{all} the minors (including $1$)
of the $r \times s$-matrix $(\nu_{\rho\sigma})_ {1\leq \rho \leq
r, 1 \leq \sigma \leq s}$; for instance, $x:= 2^{1/2}3^{1/4} +
2^{1/3} 3^{1/2}$ has degree $12$ over $\QQ$ and $\QQ[x] =
\QQ[2^{1/2}3^{1/4},2^{1/3} 3^{1/2}]$ $(= \QQ[2^{1/6}3^{1/4}]$),
cf. \cite{albu}, Example 9.2.9.

\smallskip

In a similar way, the finite subextensions of $K[\I(K_+^\times)]$
can be described for a \emph{finite} real number field $K$: The
multiplicative group $K_+^\times$ of the positive numbers in $K$
is free. (For any finite number field $K$ the group
$K^\times/\mu(K)$ ($\mu(K):=$ group of roots of unity in $K$) is
free.) If a basis $\pi_i$, $i \in I$, of $K_+^\times$ is given
(and such a basis can be constructed in principle) one has
completely analogous results for $K$ instead of $\QQ$, replacing
the primes $p \in \PP$ by the $\pi_i$, $i\in I$. (Even the
assumption that $K$ is a real field is not essential. One replaces
$K_+^\times$ by $K^\times/\mu(K)$.)

\smallskip
Iterating the construction $K\langle \{ \pm  1\} \I (K_+ ^\times)
\rangle$ from $K$, we get a tower of fields $K = K_0 \subseteq K_1
\subseteq K_2 \subseteq \ldots \subseteq \bar{K}_{\real}$ with
$K_{j+1} = K_j[\I(K_{j,+}^\times)] = K_j \langle \{\pm1\}
\I(K_{j,+}^\times) \rangle$ for an ordered field $K$. It is an
interesting task to determine for a given $x \in \bigcup_j K_j$ the
smallest $j \in \NN$ with $x \in K_j$.
\end{myexam}

\begin{myexam}\label{injhullqexam}
Any essential group extension $\QQ^\times \hookrightarrow U$ can
be embedded into the injective hull $\I(\QQ^\times) = \I(\{\pm
1\}) \times \I(\QQ_+^\times)$ and hence the universal algebra
$\QQ\langle U \rangle$ into $\QQ \langle \I(\QQ^\times) \rangle$.
We use the canonical identification $\I(\QQ^\times) = \I(\{\pm 1
\}) \times \I(\QQ_+^\times) = S^1[2^\infty] \times
\Tor(\RR_+^\times|\QQ_+^\times) \subseteq S^1 \times \RR_+^\times
= \CC^\times$, $S^1:=\{z \in \CC: |z| = 1\}$. The group
$\I(\QQ_+^\times) = \Tor(\RR_+^\times|\QQ_+^\times)$ is
torsion-free and divisible with the primes $p \in \PP$ as
canonical $\QQ$-basis and was studied in the previous example.

\smallskip
The universal algebra $\QQ\langle \I(\QQ^\times) \rangle$ is
\emph{not} a field because of $i \in S^1[2^\infty] \subseteq
\I(\QQ^\times)$, $i \notin \QQ^\times$ and $(1+i) = \zeta_8
\sqrt{2} \in \I(\QQ^\times)$, $(1+i)^4 = -4$, cf. Theorem
\ref{maintheo}.

\smallskip
To understand $\QQ\langle \I(\QQ^\times) \rangle$ we compare this
algebra with the universal $\QQ[i]$-algebra $\QQ[i] \langle
\I(\QQ[i]^\times)\rangle$ which is by Theorem \ref{maintheo} a
field.

\smallskip

Also $\I(\QQ[i]^\times)$ can be identified with a subgroup of
$\CC^\times$ which extends the identification of $\I(\QQ^\times)$
as a subgroup of $\CC^\times$ from above. We have to choose
$\tor(\I(\QQ[i]^\times)) = S^1[2^\infty]$ and take for the primes
$q \in \ZZ[i]$ with $-\pi/4 < \arg q < \pi/4$ the element
$\exp(\alpha \ln q)$ as $q^\alpha$, $\alpha \in \QQ$, and identify
$p^\alpha \in \I(\QQ^\times)$, $\alpha \in \QQ$, $p \geq 3$ prim
in $\ZZ$, in the natural way with $p^\alpha \in
\I(\QQ[i]^\times)$. For the prime $1+i \in \ZZ[i]$ and for $2 =
(-i)(1+i)^2 \in \ZZ$ we proceed as follows: $(1+i)^\alpha$,
$\alpha \in \QQ$, we will identify with $\exp(2\pi i (\alpha/8)_2)
2^{\alpha/2}$, where $r_2$ for $r \in \QQ$ denotes the
$2$-component of $[r] \in \QQ/\ZZ = \bigoplus_{p\in \PP}
(\QQ/\ZZ)[p^\infty] = \bigoplus_{p\in \PP} (\ZZ_{(p^k \!,\,k \in
\, \NN)}/\ZZ)$. Then $1+i$ will be identified with $\exp(2\pi i
/8) \sqrt{2} = 1+i$ (and hence $(1+i)^n$ with $(1+i)^n$ for all $n
\in \ZZ$). The element $2^\alpha \in \I(\QQ^\times)$, $\alpha \in
\QQ$, has in $\I(\QQ[i]^\times)$ the representation $2^\alpha =
\exp(-2\pi i (\alpha/4)_2) (1+i)^{2\alpha}$.

\smallskip

\emph{The kernel of the universal homomorphism
$\varphi:\,\QQ\langle \I(\QQ^\times) \rangle \rightarrow \QQ[i]
\langle \I(\QQ[i]^\times) \rangle$ is the principal ideal
generated by $f:= x^2 - 2x + 2 = (2i + 2)- \zeta_8 \sqrt{2}$, with
$x:= \zeta_8 \sqrt{2} \in \QQ \langle \I(\QQ^\times) \rangle =
\QQ[i] \langle \QQ[i]^\times * \I(\QQ^\times)\rangle$ and $x^4 =
-4$} (where $*$ denotes the multiplication in $\QQ \langle
\I(\QQ^\times) \rangle$ which has to be distinguished from the
multiplication in $\QQ[\I(\QQ^\times)]\subseteq \CC$). This
assertion follows from the fact that $fx_j$, $j \in J$, generate
$\ker \varphi$ as a $\QQ[i]$-vector space if $x_j$, $j \in J$,
represent the elements of the factor group $\QQ[i]^\times *
\I(\QQ^\times)/\QQ[i]^\times$. \emph{Therefore $\QQ \langle
\I(\QQ^\times) \rangle/f \QQ \langle \I(\QQ^\times) \rangle$ is
isomorphic to the subfield $\QQ[\I(\QQ^\times)] \subseteq \CC$.}
The principal ideal $(f)$ can also be generated by the idempotent
element $e:= (x+2)f/8$. If we use the automorphism of
$\I(\QQ^\times)$ induced by taking the $5$-th power on the
component $S^1[2^\infty]$ of $\I(\QQ^\times)$ and the identity on
the other components we get an automorphism $\Psi:\, \QQ \langle
\I(\QQ^\times) \rangle \rightarrow \QQ \langle \I(\QQ^\times)
\rangle$. The kernel of the homomorphism $\varphi \circ
\Psi^{-1}:\, \QQ\langle \I(\QQ^\times)\rangle \rightarrow
\QQ[i]\langle\I(\QQ[i]^\times) \rangle$ is gene\-rated by $\Psi(e)
= (-x+2)\Psi (f)/8 = 1-e$.

\smallskip

It follows that \emph{$\QQ\langle \I(\QQ^\times)\rangle$ is the
product of two fields which are both isomorphic to
$\QQ[\I(\QQ^\times)] \subseteq \CC$.} For any essential group
extension $U$ of $\QQ^\times$ we have inclusions $\QQ^\times
\subseteq U \subseteq \I(\QQ^\times)$. Hence: \emph{If $\QQ\langle U
\rangle$ is not a field, i.e. if $-4 \in U^4$ then $\QQ \langle U
\rangle$ decomposes into two fields.} But, these fields are not
necessarily isomorphic. Perhaps the simplest example is $\QQ\langle
U \rangle:= \QQ[X]/(X^{16} + 4) \cong (\QQ[X]/(X^8-2X^4+2)) \times
(\QQ[X]/(X^8 + 2X^4+2 )) = K_1 \times K_2$, $K_1 \not\cong K_2$. To
prove this, one computes for instance the Galois group $\gal(L|K)$
of the splitting field $L$ of $X^{16}+4$ over $\QQ$ and considers
$K_1$ and $K_2$ as subfields of $L$. The Galois group is isomorphic
to the semidirect product $(\ZZ_4 \times \ZZ_4) \rtimes \ZZ_2$ where
$\ZZ_2$ is generated by the complex conjugation $\kappa$ which
operates on $\ZZ_4 \times \ZZ_4$ as the matrix
\[ \left( \begin{array}{cc} 2 & 1\\ 1 & 2 \end{array} \right). \]
The two factors of the product group $\ZZ_4 \times \ZZ_4$ (which are
not conjugated in $ (\ZZ_4 \times \ZZ_4) \rtimes \ZZ_2$) are the
subgroups belonging to $K_1$ and $K_2$.
\end{myexam}

\section{Unitarily graded Galois extensions}

We consider finite Galois field extensions $L|K$. (We leave the
easy generalisations to infinite Galois extensions to the reader.
One simply uses the fact that in the graded case $L = K\langle U
\rangle = \varinjlim K\langle U^\prime \rangle$ where $U^\prime$
runs through the subgroups $U^\prime \subseteq U = {\hu {L}}$ with
$K^\times \subseteq U^\prime$ and $[U^\prime : K^\times]<
\infty$.) Let us start with the case where the Galois group is
cyclic. If $L$ has a unitary grading over $K$ then the grading
group $D$ is necessarily also cyclic. To prove this, observe that
any subgroup $D^\prime \subseteq D$ defines the graded subfield
$L_{D^\prime}$. Therefore, for any divisor $d^\prime$ of $|D| =
\ord D$, there exists at most one subgroup of $D$ of order
$d^\prime$. But this condition characterises the finite cyclic
groups in the class of all finite (not necessarily abelian) groups
$D$ (indeed, it suffices to consider prime powers $d^\prime$
dividing $|D|$). Moreover, if the cyclic extension $L|K$ has a
grading then this grading is even co-Galois and hence essentially
unique (in the sense that the group of homogeneous units is
unique). Conversely, if a Galois extension has a co-Galois grading
with cyclic grading group then the Galois group is also cyclic.
More generally, the following is true.
\ifthenelse{\boolean{optionalremark}}{(See also Remark
\ref{storchremark} for generalisations.)}{}

\begin{lem}
Let $L|K$ be a finite Galois field extension with a $D$-co-Galois grading.
Then $\exp(D) = \exp(\gal(L|K))$ and there is an element $\sigma \in \gal(L|K)$ with
$\ord \sigma = \exp(\gal(L|K))$.
\end{lem}

\begin{proof}
Let $\sigma \in G:= \gal(L|K)$. Then $L$ is graded over the
$\sigma$-invariant field $L^\sigma = L_{D^\prime}$ with grading
group $D/D^\prime$ for some subgroup $D^\prime \subseteq D$. The
extension $L|L^\sigma$ is cyclic of degree $\ord \sigma$. It follows
that $D/D^\prime$ is also cyclic of the same order. This proves
$\exp(G)|\exp(D)$. For the converse let $D^\prime \subseteq D$ be a
subgroup with cyclic quotient $D/D^\prime$ of order $\exp(D)$. Then
the Galois extension $L|L_{D^\prime}$ has a $D/D^\prime$-co-Galois
grading. By the remark above, $\gal(L|L_{D^\prime}) \subseteq
\gal(L|K) = G$ is cyclic of order $|D/D^\prime| = \exp(D)$.
\end{proof}

If the (finite) Galois extensions $L_\sigma|K$, $\sigma =
1,\ldots,s$, have a co-Galois grading and if $L:=L_1 \otimes_K
\cdots \otimes_K L_s$ is a field (i.e. if the $L_\sigma$ are
linearly disjoint over $K$), then the grading of $L$ derived from
the gradings of the factors is also co-Galois. This follows
immediately from the fact that for this grading of $L$ the
conditions of Theorem \ref{cogaloiscrit} are fulfilled since these
conditions are fulfilled for the factors. Note that a $D$-graded
Galois extension contains a root of unity of order $m$ if $D$
contains an element of order $m$. (In general, the product $L_1
\otimes_K L_2$ of co-Galois extensions is not co-Galois even if
$L_1, L_2$ are linearly disjoint. For example, $\QQ[\sqrt[3]{2}]
\otimes_\QQ \QQ[\zeta_3]$ has no co-Galois grading at all.)

\smallskip

Let us now assume that the extension $L|K$ is abelian with Galois
group $G:= \gal(L|K)$ and that it has a co-Galois grading with $U
= {\hu {L}}$ as group of homogeneous units and grading group $D
\cong U/K^\times$. Then we can prove a little bit more. If $D =
D_1 \times \cdots \times D_r$ is a decomposition of $D$ into
cyclic factors $D_\rho$, $\rho = 1,\ldots,r$, then the subfields
$L_{D_\rho}$ are also co-Galois and Galois. Hence the Galois group
$G_\rho := \gal(L_{D_\rho}|K)$ is cyclic too and $G_\rho \cong
D_\rho$. It follows the (non canonical) isomorphism $G =
\gal(L_{D_1} \otimes_K \cdots \otimes_K L_{D_r}|K)= G_1 \times
\cdots \times G_r \cong D_1 \times \cdots \times D_r = D$ (compare
also \cite[Theorem 2.9]{stefan}). Conversely, if the grading group
$D$ of an arbitrary unitary grading of an (abelian) extension
$L|K$ is isomorphic to the Galois group $G$, then the grading is
co-Galois because the mapping $D^\prime \mapsto
\gal(L|L_{D^\prime})$ is an injective and hence bijective map from
the set of subgroups $D^\prime \subseteq D$ into the set of
subgroups of $G$.

\smallskip

A (not necessarily abelian) Galois extension $L$ of $K$ which has a co-Galois
grading contains necessarily a root of unity of order $n$ where $n:= \exp(D)
= \exp(\gal(L|K))$. The base field $K$ contains necessarily a root of unity of
order $p$ for every prime divisor $p$ of $n$ and moreover a root of unity of
order $4$ if $4|n$, cf. Theorem \ref{cogaloiscrit}. Altogether, $K$
contains a root of unity of order $\ered(n)$ where $\ered(n)$ is the
\emph{extended reduction} of $n$ defined by
\[\ered(n):= \begin{cases} \red(n), \mbox{ if } 4\!\!\not| n,\\2\red(n),
\mbox{ if } 4|n.
\end{cases} \]
Here the \emph{reduction} $\red(n)$ of $n$ is the product of the prime
factors of $n$.

\smallskip

The elements of the Galois group $G$ of $L|K$ are explicitly given by the
formula
$$\sigma_\chi \left (\sum\nolimits_{d} x_d \right) = \sum\nolimits_{d} \chi(d) x_d,$$
where the index $\chi$ runs through the character group $\check D =
\Hom(D,L^\times)= \Hom(D,\mu_n(L))$, $n = \exp(D)$. \emph{It follows
that $\mu_n(L) \subseteq U = \hu {L}$ since $\sigma_\chi(U) = U$ for
all $\chi \in \check D$ $($the co-Galois grading is essentially
unique!$\,)$ and hence $\chi(d) =\sigma_\chi (x_d)/x_d \in U$ for
all $\chi \in \check D$ and all homogeneous units $x_d$ of degree
$d$, $d \in D$.}

\smallskip

The group $U = {\hu {L}}$ can be described in the following way
using only the Galois group $G$: \emph{It is}
$$U/\mu_n(L) = (L^\times/\mu_n(L))^G$$
(where $n = \exp(G)$ and the operation of $G$ on $L^\times/\mu_n(L)$ is
induced by the Galois operation). We have only to show the inclusion
$U^\prime \subseteq U$, where $U^\prime \subseteq L^\times$ is defined
by the equation $U^\prime/\mu_n(L) = (L^\times/\mu_n(L))^G$. From the
exact sequence of group cohomology
$$1 \!\!\rightarrow\! \mu_n(L)^G\! =\!\mu_n(K) \!\rightarrow\! (L\!^\times)^ G\!=\!K^\times
\!\rightarrow\!  (L\!^\times\!/\!\mu_n(L)\!)^G\!
=\!U^\prime\!/\!\mu_n(L) \!\rightarrow \!\Ho^1(G,\mu_n(L)\!)$$ we
derive the exact sequence
$$ 1 \longrightarrow \mu_n(L)/\mu_n(K) \longrightarrow
U^\prime/K^\times \longrightarrow \Ho^1(G,\mu_n(L)).$$ It follows
that $U^\prime/K^\times$is a finite group since $\Ho^1(G,\mu_n(L))$
is finite. Moreover, the exponent of $\Ho^1(G,\mu_n(L))$ divides $n
= \exp(G) = |\mu_n(L)|$.

\smallskip

\emph{We show that the universal algebra $K\langle U^\prime
\rangle$ is a field} and use Theorem \ref{maintheo} to do this. If
$p$ is a prime divisor of $|U^\prime/K^\times|$ then $p$ divides
$n = |\mu_n(L)|$ hence $\ered(n)$, and $K$ contains a root of
unity of order $p$. This proves that $K^\times \hookrightarrow
U^\prime$ is essential. If $U^\prime$ contains an element $i$ of
order $4$ but $i\notin K^\times$ then $4 \!\!\! \not | n$ (because
$|\mu_{\ered(n)}(K)|=\ered(n)$ and hence $|\mu_n(L)/\mu_n(K)|$ is
odd and $\Ho^1(G,\mu_n(L))$ does not contain an element of order
$4$. Then, by the exact sequence above, $U^\prime/K^\times$
contains no element of order $4$. It follows $-4 \notin
{U^\prime}^4$. The canonical homomorphism $K\langle U^\prime
\rangle \rightarrow L$ which extends the isomorphism $K \langle U
\rangle \stackrel{\sim} \rightarrow K[U] = L$ is injective. This
yields $U = U^\prime$.

\smallskip

We notice:

\begin{lem}\label{uprime}
Let $L|K$ be a finite Galois field extension with Galois group $G$ and
$n:= \exp(G)$. If $|\mu_n(L)| = n$, $|\mu_{\ered(n)}(K)|= \ered(n)$
and $U^\prime \subseteq L^\times$ is the subgroup with $\mu_n(L)
\subseteq U^\prime$ and $U^\prime/\mu_n(L) = (L^\times/\mu_n(L))^G$
then the universal algebra $K\langle U^\prime \rangle$ is a field
isomorphic to $K[U^\prime]\subseteq L$. The canonical
sequence
$$1\longrightarrow \mu_n(L)/\mu_n(K)\longrightarrow U^\prime/K^\times
\longrightarrow \Ho^1(G,\mu_n(L)) \longrightarrow 1$$ is exact and
$K\langle U^\prime \rangle \cong K[U^\prime]$ is a co-Galois and
Galois extension of $K$. Moreover, $K[U^\prime]\subseteq L$ is the
largest Galois subextension of $L$ which is co-Galois.
\end{lem}

\begin{proof}
The exact sequence follows from the exact sequence
$1 \rightarrow \mu_n(L) \rightarrow L^\times \rightarrow
L^\times/\mu_n(L) \rightarrow 1$ and $\Ho^1(G,L^\times) = 1$
(Noether's Theorem). The co-Galois property follows from
Theorem \ref{cogaloiscrit}. The extension $K[U^\prime]$ is
Galois since $U^\prime$ is $G$-invariant.
\end{proof}

In general, the co-Galois extension $K\langle U^\prime \rangle \cong
K[U^\prime] \subseteq L$ of Lemma \ref{uprime} is a proper subfield
of $L$. It coincides with $L$ if and only if $|U'/K^\times|=|G|$ or
equivalently
$$|\Ho^1(G,\mu_n(L))| = |\mu_n(K)|\, |G|/n.$$

This proves

\begin{theo}\label{cogaloisgaloisprop}
Let $L|K$ be a finite Galois field extension with Galois group $G$
and $n:= \exp(G)$. $L$ has a co-Galois grading over $K$ if and only
if the following conditions are satisfied:
\begin{enumerate}
\item $|\mu_n(L)| = n$ and $|\mu_{\ered(n)}(K)| = \ered(n)$. \item
$|\Ho^1(G,\mu_n(L))| = |\mu_n(K)|\, |G|/n$.
\end{enumerate}
\end{theo}

\vskip 0,2 cm

In the cyclic case condition (1) in \ref{cogaloisgaloisprop} is sufficient:

\begin{theo}\label{cycliccase}
Let $L|K$ be a finite cyclic field extension of degree $n$. Then
$L$ has a unitary grading $($which is necessarily a co-Galois
grading$)$ if and only if $|\mu_n(L)| = n$ and
$|\mu_{\ered(n)}(K)| =\ered(n)$.
\end{theo}

\begin{proof}
Let the conditions on the roots of unity be satisfied. We have to
prove that condition (2) of Theorem \ref{cogaloisgaloisprop} is also
satisfied, which means $|\Ho^1(\gal(L|K),\mu_n(L))| =|\mu_n(K)|$.
Let $\sigma$ be a generator of the Galois group $G:=\gal(L|K)$. Then
the cohomology group $\Ho^1(G,\mu_n(L))$ is the homology of the
complex
$$\mu_n(L) \stackrel{\sigma/\id}\longrightarrow \mu_n(L)
\stackrel{\norm}\longrightarrow \mu_n(L)$$ of finite groups where
$\norm$ is the norm $x \mapsto \prod_{j=0}^{n-1} \sigma^j x$.
It follows from the index satz
$$|\Ho^1(G,\mu_n(L))| = |\ker \sigma/\id|\, |\coker \norm|/|\mu_n(L)|
= |\mu_n(K)|\, |\coker \norm| /n.$$ It remains to show: $|\coker
\norm| = n$, i.e. $\mu_n(L)$ belongs to the norm-1-group of $L|K$.
But this is verified by the following (probably well-known) lemma.
\end{proof}

\begin{lem}\label{rootnorm}
Let $L|K$ be a finite field extension of degree $n$. Then
$\mu_n(L)$ is contained in the norm-$1$-group of $L|K$.
\end{lem}

\begin{proof}
It is sufficient to show: If $\zeta \in L$ is a root of unity of
prime power order $p^\alpha>1$ and if $p^\alpha$ divides $n$, then
$\norm^L_K (\zeta) = 1$. Consider the subfield $K[\zeta]$ and let
$m:=[K[\zeta]:K]$. Then $m|n$ and $\norm^L_K(\zeta) =
\norm^{K[\zeta]}_K(\norm^L_{K[\zeta]} (\zeta)) =
\norm^{K[\zeta]}_K(\zeta^{n/m})$ and $\zeta^{n/m} \in \mu_m
(K[\zeta])$. Therefore, we may assume additionally $L=K[\zeta]$.
$K[\zeta]|K$ is a Galois extension. Its Galois group is a subgroup
of the automorphism group $\aut(\langle \zeta \rangle) = (\ZZ/\ZZ
p^\alpha)^\times$ and its order $m$ divides $p^{\alpha -1}(p-1)$,
i.e. $m = p^\beta t$, $\beta< \alpha$, $t|(p-1)$.

\smallskip

It is enough to prove $\norm(\zeta)^{p^{\alpha-\beta}} :=
\norm^{K[\zeta]}_{K[\zeta^{p^{\alpha-1}}]} (\zeta)^{p^{\alpha-\beta}} = 1$.
$K[\zeta]|K[\zeta^{p^{\alpha-1}}]$ is a Galois extension of degree
$p^\beta$ and its Galois group $G$ is a subgroup of $1 + \ZZ p/ \ZZ p^\alpha
\subseteq (\ZZ / \ZZ p^\alpha)^\times$.

\smallskip

First let $p \geq 3$. Then $G = 1 + \mathfrak a$,
$\mathfrak a:=\ZZ p^{\alpha-\beta}/ \ZZ p^\alpha$ and $\norm(\zeta)
^{p^{\alpha-\beta}} = \left( \prod_{\sigma \in G} \sigma \zeta\right)
^{p^{\alpha-\beta}} = \zeta^{p^{\alpha-\beta}S}$, $S:= \sum_{j \in \mathfrak a}
(1+j) = p^\beta + \sum_{j \in \mathfrak a} j = p^\beta$ since
$\sum_{j \in \mathfrak a} j = 0$, hence $\norm(\zeta)^{p^{\alpha-\beta}}
= \zeta^{p^{\alpha-\beta}p^\beta} = 1$.

\smallskip

Now let $p = 2$ and $\alpha \geq 2$. Then $1 + \ZZ 2/ \ZZ 2^\alpha$ is the
product of the cyclic subgroups $\{\pm 1\}$ and $1 + \ZZ 4/ \ZZ 2^\alpha$.
The subgroups of order $2^\beta$ are $1+ \mathfrak a$, $\mathfrak a :=
\ZZ 2^{\alpha-\beta} / \ZZ 2^\alpha$ (if $\beta \leq \alpha-2$) and the groups
$(1+ \alpha) \uplus -(1+ \mathfrak a) (1+x)$ with $\mathfrak a := \ZZ 2
^{\alpha-\beta+1} / \ZZ 2^\alpha$ and a fixed $x \in \ZZ 4 / \ZZ 2^\alpha$,
$(1+x)^2 -1 = x(2+x) \in \mathfrak a$.

\smallskip

In the first case $\norm(\zeta)^{2^{\alpha-\beta}} = \zeta^{2^{\alpha-\beta}S}$
with $S:= \sum_{j \in \mathfrak a} (1+j) = 2^\beta + \sum_{j \in \mathfrak a}
j = 2^\beta + 2^{\alpha-1}$ if $\beta>0$ (and $S=1$ if $\beta =0$), hence
$\norm(\zeta)^{2^{\alpha-\beta}} = 1$. In the second case $\norm (\zeta)^{2^{\alpha-\beta}} = \zeta
^{2^{\alpha-\beta}S}$ with $S:= -(\sum_{j\in \mathfrak a} j) x - 2^{\beta-1}x$,
hence $\zeta^{2^{\alpha-\beta}S} = \zeta^{-2^{\alpha-1}x} = 1$.
\end{proof}

In general, condition (1) in Theorem \ref{cogaloisgaloisprop} is
not sufficient for the existence of a co-Galois grading of $L|K$,
even in the abelian case. For instance, the Galois extension
$\QQ[\sqrt{-3},\sqrt{-19}] \subseteq \QQ[\zeta_{3^2\cdot 19}]$
with Galois group $\ZZ_3 \times \ZZ_9$ has no co-Galois grading
but $\zeta_9 \in \QQ[\zeta_{3^2\cdot 19}]$ and $\zeta_3 \in
\QQ[\sqrt{-3},\sqrt{-19}]$. -- If the Galois group $G$ of $L|K$ is
abelian and contains a subgroup isomorphic to $\ZZ_n \times
\ZZ_n$, $n =\exp(G)$, then $L|K$ has a co-Galois grading (if and)
only if $L|K$ is a Kummer extension, i.e. $|\mu_n(K)| = n$.\,--
Also, if $|\mu_n(K)|=n$ and $L|K$ has a co-Galois grading then $G$
is necessarily abelian, hence $L|K$ is a Kummer extension. It
follows, quite generally, that for a finite Galois and co-Galois
extension $L|K$ with Galois group $G$ the co-Galois extension
$L|K[\zeta_n]$ ($n=\exp(G)$) is a Kummer extension. Since an
abelian extension $L|K$ has a co-Galois grading if and only if
every cyclic subextension $L^\prime|K$, $L^\prime \subseteq L$,
has such a grading Theorem \ref{cycliccase} is useful also in this
more general setting.

\smallskip

With respect to Lemma \ref{rootnorm} the group $\mu_n(K) =
\Ho^0(G,\mu_n(L))$ can also be interpreted as the modified
cohomology group $\hat \Ho^0(G, \mu_n(L))$ (in the sense of Tate).
Since $\hat \Ho^1 (G,\mu_n(L)) = \Ho^1(G,\mu_n(L))$ condition (2) in
\ref{cogaloisgaloisprop} can be written as
$$\her(G,\mu_n(L)):= \frac{|\hat \Ho^0(G,\mu_n(L))|}{|\hat \Ho^1(G,\mu_n(L))|}
=\frac{n}{|G|}\,,$$ $n = \exp(G)$. Since for $G$ cyclic and for the
finite $G$-module $\mu_n(L)$, the quotient $\her(G,\mu_n(L))$ (which
is called the Herbrand quotient) is always $1$, we get Theorem
\ref{cycliccase} in a more conceptual way. Let us also mention the
classical description of the cohomology group $\hat \Ho^1
(G,\mu_n(L)) = \Ho^1(G,\mu_n(L))$ as
$$\hat \Ho^1 (G,\mu_n(L)) = {L^\times}^n \cap K^\times/{K^\times}^n$$
derived from the exact sequence $1 \rightarrow \mu_n(L) \rightarrow
L^\times \stackrel{n} \rightarrow {L^\times}^n \rightarrow 1$ and
$\hat \Ho^1(G, L^\times) = 1$.

\smallskip

If $L|K$ is an extension of \emph{finite} fields with $|K| = q$ and
$|L| = q^n$ then $|\mu_{\ered(n)}(K)| = \ered(n)$ is equivalent with
$q \equiv 1 \!\!\mod\! \ered(n)$. Of course, this condition implies
$q^n \equiv 1 \!\! \mod \!n$, i.e. $|\mu_n(L)|=n$. Theorem
\ref{cycliccase} has therefore the following corollary which can
also be proved more directly.

\begin{cor}
An extension $L|K$ of finite fields of degree $n$ with $q = |K|$ has
a unitary grading if and only if $q \equiv 1\!\! \mod \! \ered(n)$.
In this case, the grading is a co-Galois grading with cyclic grading
group and in particular essentially unique.
\end{cor}

\vskip 0,2 cm

\begin{myexam}
A cyclotomic field $\QQ[\zeta_n]$ over $\QQ$ can have a co-Galois
grading only in the case $\ered(\varphi(n)) \leq 2$ which implies
$n | 24$. In this case it has a co-Galois grading for trivial
reasons (cf. also \cite[Corollary 7.4.5]{albu}).

\smallskip

A little bit more complicated is the problem to determine the $n$
for which $\QQ[\zeta_n]|\QQ$ has a unitary (not necessarily
co-Galois) grading. \emph{This is the case exactly for}
$$n = 2^\alpha \cdot 3^\beta,~\alpha \in \NN,~\beta \in \{0,1\}.$$
To see this, one can use the following strategy (for a more detailed
account see \cite{kaiddiplom}): Let $L:=\QQ[\zeta_n]$ be a
cyclotomic field which is unitarily $D$-graded over $\QQ$. First
consider the case that $n = p^\alpha$ is a prime power. For $p = 2$
there is nothing to prove, so let $p \geq 3$. By considering roots
of unity one gets $(p^\alpha - p^{\alpha - 1})|2p^\alpha$ which
yields $p = 3$. Because the cyclic extension $\QQ[\zeta_9]|\QQ$
contains the real subfield $\QQ[\zeta_9] \cap \RR$ of degree $3$
over $\QQ$ we get $n = 3$.

\smallskip

Now we treat the general case. We can assume that $n$ is even, $n
> 2$ and $\varphi(n)| n$. We show that $\varphi(n)$ has to be a
power of $2$, i.e. $n = 2^\alpha p_1 \cdots p_r$ with Fermat primes
$p_j$, $j=1,\ldots,r$. Assume there is an odd prime divisor $p$ of
$\varphi(n)$. Then there exists a subgroup $\tilde{D}$ of $D$ of
order $p$ and $\QQ[\zeta_p] \subseteq L_{\tilde{D}}$. But this is a
contradiction. Hence the grading group $D$ is a $2$-group and
moreover $\exp(D) \leq 2^\alpha$. Now let $D = D_1 \times \cdots
\times D_s$ be a decomposition of $D$ into cyclic groups. Then the
subfields $L_{D_j}$, $j = 1,\ldots,s$, are linearly disjoint over
$\QQ$. Hence $D$ has to be of the form $D \cong \ZZ_{2^e} \times
\ZZ_2 \times \cdots \times \ZZ_2$ with $2^e = \exp(D)$. This yields
also $\exp(\gal(L|\QQ)) \leq \exp(D)$.

\smallskip

If $\alpha = 1$ we get obviously $n = 6$ and for
$\alpha = 2$ one easily checks that $n = 4$ or $n = 12$. Now let
$\alpha \geq 3$. By comparing the Galois group
$$\gal(L|\QQ) = (\ZZ / \ZZ n)^\times \cong \ZZ_2 \times \ZZ_{2^{\alpha-2}}
\times \ZZ_{p_1 - 1} \times \cdots \times \ZZ_{p_r - 1}$$ and the
grading group $D$ one finds that $n = 2^\alpha (\cdot 3) \cdot 5$
(the factor $3$ is optional) and $\exp(D) = 2^\alpha$ is the only
critical case. Then we consider the tower of fields $\QQ \subseteq
\QQ[\zeta_{2^\alpha}] \subseteq L_{\ZZ_{2^\alpha}} \subseteq L$. By
Galois theory we see that $L_{\ZZ_{2^\alpha}} \cap \QQ[\zeta_5] =
\QQ[\sqrt{5}]$. Hence $E:= \QQ[i,\sqrt{2},\sqrt{5}] \subseteq
L_{\ZZ_{2^\alpha}}$ and $\gal(E|\QQ) \cong \ZZ_2 \times \ZZ_2 \times
\ZZ_2$. But this is a contradiction to $\gal(L_{\ZZ_{2^\alpha}}|\QQ)
\cong \ZZ_{2^{\alpha-2}} \times \ZZ_4$.
\end{myexam}

\bibliographystyle{plain}

\begin{thebibliography}{10}

\bibitem{albu}
T.~Albu.
\newblock {\em Cogalois Theory}.
\newblock New\! York Basel, 2002.

\bibitem{besicovitch}
A.~Besicovitch.
\newblock On the linear independence of fractional powers of integers.
\newblock {\em J. London Math. Soc.}, 15:3--6, 1940.

\bibitem{dade}
E.~C. Dade.
\newblock Group-graded rings and modules.
\newblock {\em Math. Z.}, 174:241--262, 1980.

\bibitem{greitherharrison}
C.~Greither and D.~K. Harrison.
\newblock A {Galois} correspondence for radical extensions of fields.
\newblock {\em J. Pure Appl. Algebra}, 43:257--270, 1986.

\bibitem{kaiddiplom}
A.~Kaid.
\newblock {Unit\"ar-graduierte K\"orpererweiterungen}.
\newblock Diplomarbeit Bo\-chum, 2004.

\bibitem{kneserlinearwurzel}
M.~Kneser.
\newblock {Lineare Abh\"angigkeit von Wurzeln}.
\newblock {\em Acta Arith.}, 26:307--308, 1975.

\bibitem{schejastorch2}
G.~Scheja and U.~Storch.
\newblock {\em Lehrbuch der Algebra, Teil 2}.
\newblock Teubner, 1988.

\bibitem{schinzel}
A.~Schinzel.
\newblock On linear dependence of roots.
\newblock {\em Acta Arith.}, 28:161--175, 1975.

\bibitem{siegelwurzel}
C.~L. Siegel.
\newblock {Algebraische Abh\"angigkeit von Wurzeln}.
\newblock {\em Acta Arith.}, 21:59--64, 1972.

\bibitem{stefan}
D.~Stefan.
\newblock Cogalois extensions via strongly graded fields.
\newblock {\em Comm. Algebra}, 27:5687--5702, 1999.

\end{thebibliography}

\noindent Holger Brenner (h.brenner@shef.ac.uk), Almar Kaid
(a.kaid@shef.ac.uk)

\noindent Department of Pure Mathematics, University of Sheffield,
Hicks Building, Hounsfield Road, Sheffield S3 7RH, United Kingdom.
\medskip

\noindent Uwe Storch (Uwe.Storch@ruhr-uni-bochum.de)

\noindent
 Fakult\"at f\"ur Mathematik der Ruhr-Universit\"at
Bochum, Universit\"atsstra\ss e 150, D-44801 Bochum, Germany.

\end{document}